\documentclass{amsart}
\usepackage[dvips]{epsfig}
\usepackage{graphicx}
\usepackage{latexsym}
\usepackage{amsmath}
\usepackage{amsthm}
\usepackage{amssymb}
\usepackage{mathrsfs}
\usepackage[shortlabels]{enumitem}
\newlist{propenum}{enumerate}{1} 
\setlist[propenum]{label=(\roman*)}

\usepackage{eucal}

\usepackage{xcolor}

\usepackage{mathtools}

 \usepackage[applemac]{inputenc}

\newtheorem{thm}{Theorem}[section]
\newtheorem{theo}[thm]{Theorem}

\newtheorem{lem}[thm]{Lemma}

\newtheorem{defi}[thm]{Definition}
\newtheorem{hyp}[thm]{Assumption}

\theoremstyle{remark}

\newtheorem{rem}[thm]{Remark}

\newcommand{\inL} {\underset{n \rightarrow \infty}{\overset{\rm (d)}{\xrightarrow{\hspace*{0.75cm}}}} }

\newcommand{\nlim} {\underset{n \rightarrow \infty}{{\xrightarrow{\hspace*{0.5cm}}}} }

\newcommand{\vip}{\vskip.2cm}
\newcommand{\COMMENTAIRE}[1]{}

\newcommand{\field}[1]{\mathbb{#1}}
\newcommand{\EE}{\field{E}}

\newcommand{\GG}{\field{G}}

\newcommand{\NN}{\field{N}}

\newcommand{\TT}{\field{T}}

\newcommand{\Bb}{{\mathcal B}}

\newcommand{\Ff}{{\mathcal F}}

\newcommand{\Pp}{{\mathcal P}}

\newcommand{\Qq}{{\mathcal Q}}
\newcommand{\vt}{{\vartriangle}}

\newcommand{\rd}{{\rm d}}

\newcommand{\bF}{{\mathfrak f}}

\newcommand{\cb}{{\mathcal B}}
\newcommand{\cc}{{\mathcal C}}

\newcommand{\ce}{{\mathcal E}}

\newcommand{\cf}{{\mathcal F}}

\newcommand{\cp}{{\mathcal P}}
\newcommand{\cq}{{\mathcal Q}}
\newcommand{\crr}{{\mathcal R}}

\newcommand{\cs}{{\mathscr S}}

\newcommand{\C}{{\mathbb C}}

\newcommand{\E}{{\mathbb E}}

\newcommand{\G}{\mathbb{G}}
\newcommand{\N}{{\mathbb N}}
\renewcommand{\P}{{\mathbb P}}

\newcommand{\R}{{\mathbb R}}

\newcommand{\T}{\mathbb{T}}

\newcommand{\sot}{\otimes_{\rm sym}}
\newcommand{\ssub}{\Sigma^{\rm sub}}
\newcommand{\stsub}{\Sigma^{\vt, {\rm sub}}}

\newcommand{\stcrit}{\Sigma^{\vt, {\rm crit}}}

\newcommand{\norm}[1]{\mathop{\parallel\! #1 \! \parallel}\nolimits}

\newcommand{\inv}[1]{\mathop{\frac{1}{ #1}}\nolimits}

\newcommand{\reff}[1]{(\ref{#1})}

\begin{document}

\title[clt for bmc:extension to the mother-daughters triangle]{Central limit theorem for bifurcating Markov chains: the mother-daughters triangles case}

\author{S. Val\`ere Bitseki Penda}

\address{S. Val\`ere Bitseki Penda, IMB, CNRS-UMR 5584, Universit\'e Bourgogne Franche-Comt\'e, 9 avenue Alain Savary, 21078 Dijon Cedex, France.}

\email{simeon-valere.bitseki-penda@u-bourgogne.fr}

\date{\today}

\begin{abstract}
The main objective of this article is to establish a central limit theorem for additive three-variable functionals of bifurcating Markov chains. We thus extend the central limit theorem under point-wise ergodic conditions studied in Bitseki-Delmas (2022) and to a lesser extent, the results of Bitseki-Delmas (2022) on central limit theorem under $L^{2}$ ergodic conditions. Our results also extend and complement those of Guyon (2007) and Delmas and Marsalle (2010). In particular, when the ergodic rate of convergence is greater than $1/\sqrt{2}$, we have, for certain class of functions,  that the asymptotic variance is non-zero at a speed faster than the usual central limit theorem studied by Guyon and Delmas-Marsalle.
\end{abstract}

\maketitle






\section{Introduction}

This article is devoted to the extension of the central limit theorem for bifurcating Markov chains given in Bitseki-Delmas \cite{BD1} when the functions do not only depend on one variable $x$ say, but on the  triplet, $xx_{0}x_{1}:=(x,x_{0},x_{1})$ say. Using the terminology given in \cite{Guyon}, we will talk about the mother-daughters triangles case. The study of mother-daughters triangles for bifurcating Markov chains models is particularly important to make statistical inference (see for e.g \cite{Guyon, DM10, bitsekidjellout2014, BO17, bdsgp:aa}). The results given here allow us to extend those obtained by Guyon \cite{Guyon} and Delmas and Marsalle \cite{DM10}. Indeed, in their works, the authors in \cite{DM10, Guyon} considered conditionally centered additive functionals (see Remark \ref{rem:ext-GDM} for more details), which is equivalent to study increments of martingale.  As in \cite{BD1}, we prove the existence of three regimes for the central limit theorem. However, these three regimes disappear when considering conditionally centered additive functionals (see Section \ref{sec:special} and in particular Theorem \ref{theo:triangPf} for more details). We stress that with appropriate hypothesis, the results obtained in this paper also hold in the $L^{2}$ ergodic case 
(see \cite{BD2} for more details).
Now, before giving the plan of the paper, we introduce useful definitions, notations and assumptions.

\subsection{Bifurcating Markov chains}
If $(E,  \ce)$ is a  measurable space, then $\cb(E)$  (resp. $\cb_b(E)$,
resp.    $\cb_+(E)$)  denotes   the  set   of  (resp.    bounded,  resp.
non-negative)  $\R$-valued measurable  functions  defined  on $E$.   For
$f\in \cb(E)$, we set $\norm{f}_\infty =\sup\{|f(x)|, \, x\in E\}$.  For
a finite measure $\lambda$ on $(E,\ce)$ and $f\in \cb(E)$ we shall write
$\langle \lambda,f  \rangle$ for  $\int f(x) \,  \rd\lambda(x)$ whenever
this  integral is  well defined.   If $(E,d)$  is a  metric space,  then
$\ce $ will denote its Borel $\sigma$-field and the set $\cc_b(E)$ (resp.
$\cc_+(E)$) denotes the set of bounded (resp.  non-negative) $\R$-valued
continuous functions defined on $E$.

Let  $(S, \cs)$   be  a  measurable  space. 
Let $Q$ be a   
probability kernel   on $S \times \cs$, that is:
$Q(\cdot  , A)$  is measurable  for all  $A\in \cs$,  and $Q(x,
\cdot)$ is  a probability measure on $(S,\cs)$ for all $x \in
S$. For any $f\in \cb_b(S)$,   we set for $x\in S$:
\begin{equation}
   \label{eq:Qf}
(Qf)(x)=\int_{S} f(y)\; Q(x,\rd y).
\end{equation}
We define $(Qf)$, or simply $Qf$, for $f\in \cb(S)$ as soon as the
integral \reff{eq:Qf} is well defined, and we have $\cq f\in \cb(S)$. For $n\in \N$, we denote by $Q^n$  the
$n$-th iterate of $Q$ defined by $Q^0=I_d$, the identity map on
$\cb(S)$, and $Q^{n+1}f=Q^n(Qf)$ for $f\in \cb_b(S)$.  

Let $P$ be a   
probability kernel   on $S \times \cs^{\otimes 2}$, that is:
$P(\cdot  , A)$  is measurable  for all  $A\in \cs^{\otimes 2}$,  and $P(x,
\cdot)$ is  a probability measure on $(S^2,\cs^{\otimes 2})$ for all $x \in
S$. For any $g\in \cb_b(S^3)$ and $h\in \cb_b(S^2)$,   we set for $x\in S$:
\begin{equation}
   \label{eq:Pg}
(Pg)(x)=\int_{S^2} g(x,y,z)\; P(x,\rd y,\rd z)
\quad\text{and}\quad
(Ph)(x)=\int_{S^2} h(y,z)\; P(x,\rd y,\rd z).
\end{equation}
We define $(Pg)$ (resp. $(Ph)$), or simply $Pg$ for $g\in \cb(S^3)$
(resp. $Ph$ for $h\in \cb(S^2)$), as soon as the corresponding 
integral \reff{eq:Pg} is well defined, and we have  that $Pg$ and
$Ph$ belong to $\cb(S)$.
\medskip 

We  now  introduce   some  notations  related  to   the  regular  binary
tree.   Recall    that  $\N$ is the set of non-negative integers and
$\N^*=    \N   \setminus   \{0\}$. We set 
$\T_0=\G_0=\{\emptyset\}$,              $\G_k=\{0,1\}^k$             and
$\T_k = \bigcup  _{0   \leq  r  \leq  k}  \G_r$  for   $k\in  \N^*$,  and
$\T  = \bigcup  _{r\in \N}  \G_r$.  The  set $\G_k$  corresponds to  the
$k$-th generation,  $\T_k$ to the tree  up the $k$-th  generation, and $\T$ the
complete binary tree. For $i\in \T$, we denote by $|i|$  the generation of $i$ ($|i|=k$ if and
only if $i\in \G_k$) and $iA=\{ij; j\in A\}$ for $A\subset \T$, where
$ij$ is the concatenation of the two sequences $i,j\in \T$, with the
convention that $\emptyset i=i\emptyset=i$.

We recall the definition of bifurcating Markov chain (BMC) from Guyon
\cite{Guyon}. 
\begin{defi}
  We say  a stochastic process indexed  by $\T$, $X=(X_i,  i\in \T)$, is
  a bifurcating Markov chain on a measurable space $(S, \cs)$ with
  initial probability distribution  $\nu$ on $(S, \cs)$ and probability
  kernel $\cp$ on $S\times \cs^{\otimes 2}$, a  BMC in
  short, if:
\begin{itemize}
\item[-] (Initial  distribution.) The  random variable  $X_\emptyset$ is
  distributed as $\nu$.
   \item[-] (Branching Markov property.) For  a sequence   $(g_i, i\in
     \T)$ of functions belonging to $\cb_b(S^3)$, we have for all $k\geq 0$,
\[
\E\Big[\prod_{i\in \G_k} g_i(X_i,X_{i0},X_{i1}) |\sigma(X_j; j\in \T_k)\Big] 
=\prod_{i\in \G_k} \cp g_i(X_{i}).
\]
\end{itemize}
\end{defi}

We    define     three
probability kernels $P_0, P_1$ and $\cq$ on $S\times \cs$ by:
\[
P_0(x,A)=\cp(x, A\times S), \quad
P_1(x,A)=\cp(x, S\times A) \quad\text{for
$(x,A)\in S\times  \cs$, and}\quad
\cq=\inv{2}(P_0+P_1).
\] 
Notice  that  $P_0$ (resp.   $P_1$)  is  the  restriction of  the  first
(resp. second) marginal of $\cp$ to $S$.  Following Guyon \cite{Guyon}, we
introduce an  auxiliary Markov  chain $Y=(Y_n, n\in  \N) $  on $(S,\cs)$
with  $Y_0$ distributed  as $X_\emptyset$  and transition  kernel $\cq$.
The  distribution of  $Y_n$ corresponds  to the  distribution of  $X_I$,
where $I$  is chosen independently from  $X$ and uniformly at  random in
generation  $\G_n$.    We  shall   write  $\E_x$   when  $X_\emptyset=x$
(\textit{i.e.}  the initial  distribution  $\nu$ is  the  Dirac mass  at
$x\in S$).  

\medskip

For all $u\in \TT$, we denote by $X_{u}^{\vartriangle} = (X_{u},X_{u0},X_{u1})$ the mother-daughters triangle. One can check that $(X_{u}^{\vartriangle},u\in\TT)$ is a bifurcating Markov chain on $S^{3}$ with transition probability $\Pp^{\vartriangle}$ defined by
\[\Pp^{\vartriangle}(xx_{0}x_{1}, dyy_{0}y_{1}, dzz_{0}z_{1}) = \delta_{x_{0}}(dy)\Pp(y,dy_{0},dy_{1})\delta_{x_{1}}(dz)\Pp(z,dz_{0},dz_{1}),\]
where we set  $xx_{0}x_{1} = (x,x_{0},x_{1}).$ The first and the second marginal of $\Pp^{\vartriangle}$ are defined by:
\begin{equation*}
P_{0}^{\vartriangle}(xx_{0}x_{1}, dyy_{0}y_{1}) = \delta_{x_{0}}(dy)\Pp(y,dy_{0},dy_{1}) \quad \text{and} \quad P_{1}^{\vartriangle}(xx_{0}x_{1}, dzz_{0}z_{1}) = \delta_{x_{1}}(dy)\Pp(z,dz_{0},dz_{1}).
\end{equation*}
We consider the sequence $(Y_{n}^{\vartriangle}, n\in \NN)$ defined recursively as follow: we set $Y_{0}^{\vartriangle} = X_{\emptyset}^{\vartriangle}$ and for all $n\in \NN$, if we are on a vertex $u \in \GG_{n}$ of the tree, then with probability $1/2$, we choose the vertex $u0$ and we set $Y_{n+1}^{\vartriangle} = X_{u0}^{\vartriangle}$ or we choose the vertex $u1$ and we set $Y_{n+1}^{\vartriangle} = X_{u1}^{\vartriangle}$. One can easily see that the sequence $(Y_{n}^{\vartriangle}, n\in \NN)$ is a Markov chain on $S^{3}$ whose transition probability $\Qq^{\vartriangle}$ is defined by
\[\Qq^{\vartriangle}(xx_{0}x_{1}, dyy_{0}y_{1}) = \frac{1}{2}(\delta_{x_{0}}(dy) + \delta_{x_{1}}(dy))\Pp(y,dy_{0},dy_{1}).\]

\medskip

Let $i,j\in \T$. We write $i\preccurlyeq  j$ if $j\in i\T$. We denote by $i\wedge j$  the most recent  common ancestor of  $i$ and $j$,  which is defined  as   the  only   $u\in  \T$   such  that   if  $v\in   \T$  and $ v\preccurlyeq i$, $v \preccurlyeq j$  then $v \preccurlyeq u$. We also define the lexicographic order $i\leq j$ if either $i \preccurlyeq j$ or $v0  \preccurlyeq i$  and $v1  \preccurlyeq j$  for $v=i\wedge  j$.  Let $X=(X_i, i\in  \T)$ be  a $BMC$  with kernel  $\cp$ and  initial measure $\nu$. For $i\in \T$, we define the $\sigma$-field:
\begin{equation*}\label{eq:field-Fi}
\cf_{i}=\{X_u; u\in \T \text{ such that  $u\leq i$}\}.
\end{equation*}
By construction,  the $\sigma$-fields $(\cf_{i}; \, i\in \T)$ are nested as $\cf_{i}\subset \cf_{j} $ for $i\leq  j$.

\medskip

By convention, for  $f,g\in  \cb(S)$,   we  define  the  function   $f\otimes  g$  by $(f\otimes g)(x,y)=f(x)g(y)$ for  $x,y\in S$ and 
\[
f\sot g= \inv{2}(f\otimes g + g\otimes f) \quad\text{and}\quad f\otimes ^2= f\otimes f.
\]

For $f \in \Bb(S^{3})$ and a finite subset $A \subset \TT$, we define:
\begin{equation*}
M_{A}(f) = \sum_{u\in A}f(X_{u}^{\vartriangle}).
\end{equation*}
In the sequel we will also use the following notation: let $g$ and $h$ be two functions which depend on one variable, $x$ say; we denote by $g\oplus h$ the function of three variables, $xx_{0}x_{1} =(x,x_{0},x_{1})$ say, defined by
\begin{equation*} 
(g\oplus h)(x,x_{0},x_{1}) = g(x_{0}) + h(x_{1}).
\end{equation*}
 
\subsection{Assumptions}
Let  $X=(X_u, u\in  \T)$ be  a BMC  on $(S,  \mathcal{S})$ with  initial probability distribution $\nu$, and probability kernel $\cp$.  Recall that    $\cq$ is the induced Markov kernel. We will present the results according to the following hypothesis. 
   
\medskip

For  a  set  $F\subset  \cb(S)$   of  $\R$-valued  functions,  we  write
$F^2=\{f^2; f\in F\} $, $F\otimes  F=\{f_0\otimes f_1; f_0, f_1\in F\}$,
and  $P(E)=\{Pf;  f\in E\}$  whenever  a  kernel $P$  act  on  a set  of
functions $E$.  We state  first a  structural assumption  on the  set of
functions we shall consider.

\begin{hyp}
   \label{hyp:F}
Let $F\subset \cb(S)$ be a set of $\R$-valued functions such that:
\begin{itemize}
   \item[$(i)$] $F$ is a  vector subspace which contains the constants;
   \item[$(ii)$] $F^2 \subset F$;
   \item[$(iii)$] $F\subset L^1(\nu)$; 
 \item[$(iv)$]  $F\otimes F \subset L^1(\cp(x, \cdot))$ for all $x\in S$,
    and $\cp(F\otimes F)\subset F$.
\end{itemize}
\end{hyp}

The   condition   $(iv)$   implies   that   $P_0(F)\subset   F$,
$P_1(F)\subset F$ as  well as $\cq(F)\subset F$.  Notice that if  $f\in
F$, then even
if $|f|$ does  not belong to $F$, using  conditions  $(i)$  and $(ii)$,
we get, with $g=(1+f^2)/2$, that   $|f|\leq  g$ and $g\in F$.
\medskip

We  consider the following  ergodic properties for $\cq$. 
\begin{hyp}
   \label{hyp:F1}
There exists a probability measure $\mu$ on $(S, \cs)$ 
  such that $F\subset L^1(\mu)$ and for all $f\in F$, we have the
  point-wise convergence  $\lim_{n\rightarrow \infty } \cq^{n}f =
  \langle \mu, f \rangle$ and 
there exists $g\in F$ with:
\begin{equation}
   \label{eq:erg-bd}
 |\cq^n(f)|\leq  g\quad\text{for all $n\in
  \N$.}
 \end{equation} 
\end{hyp}

We consider also the following geometrical ergodicity. 
\begin{hyp}
   \label{hyp:F2}
There exists a probability measure $\mu$ on $(S, \cs)$  such that
$F\subset L^1(\mu)$,  and $\alpha\in (0, 1)$ such that for all $f\in F$
  there exists $g\in F$ such that: 
 \begin{equation}
   \label{eq:geom-erg}
|\cq^{n}f - \langle \mu, f \rangle| \leq \alpha^{n} g \quad
\text{for all  $n\in \N$.}
\end{equation}
\end{hyp}

A  sequence $\bF=(f_\ell,  \ell\in  \N)$ of  elements  of $F$  satisfies
uniformly \reff{eq:erg-bd}  and \reff{eq:geom-erg} if there  is $g\in F$
such that:
\begin{equation}\label{eq:unif-f}
|\cq^n(f_\ell)|\leq  g \quad\text{and}\quad |\cq^{n}f_\ell  - \langle \mu, f_\ell \rangle|  \leq \alpha^{n} g \quad \text{for all $n,\ell\in \N$.}
\end{equation}
This    implies    in    particular   that    $|f_\ell|\leq    g$    and
$|\langle \mu, f_\ell \rangle|\leq \langle \mu, g \rangle$.  Notice that
\reff{eq:unif-f} trivially holds if $\bF$ takes finitely distinct values
(\textit{i.e.} the subset $\{f_\ell; \ell\in  \N\}$ of $F$ is finite)
each satisfying \reff{eq:erg-bd}  and \reff{eq:geom-erg}. 
\medskip

We consider the stronger ergodic property based on a second spectral
gap. 
\begin{hyp}\label{hyp:F3}
There exists a probability measure $\mu$ on $(S, \cs)$  such that $F\subset L^1(\mu)$,  and $\alpha\in (0, 1)$, a finite  non-empty set $J$ of indices, distinct complex eigenvalues $\{\alpha_j, \, j\in J\}$ of  the  operator  $\cq$ with  $|\alpha_j|=\alpha$,  non-zero  complex projectors $\{\crr_j, \, j\in J\}$  defined on $\C F$, the $\C$-vector space    spanned by $F$, such that $\crr_j\circ \crr_{j'}= \crr_{j'}\circ  \crr_{j} = 0$ for all  $j\neq j'$ (so that  $\sum_{j\in J} \crr_j$ is  also a projector defined  on $\C F$) and a positive  sequence $(\beta_n, n\in \N)$, non-increasing, bounded from above by 1 and converging  to $0$, such that  for  all  $f\in  F$  there  exists  $g\in  F$  and,  with $\theta_j=\alpha_j/\alpha$:
\[\Big|\cq^{n}( f) - \langle \mu, f \rangle - \alpha^n \sum_{j\in J}\theta_j^n\,  \crr_j (f) \Big| \leq \beta_n \alpha^{n} g \quad \text{for all  $n\in \N$.}\]
\end{hyp}

We shall  consider sequences $\bF=(f_\ell,  \ell\in \N)$ of  elements of $F$ which satisfies Assumption \ref{hyp:F3} uniformly, that is such that there exists $g\in F$ with:
\begin{equation}\label{eq:unif-f-crit}
|\cq^n(f_\ell)|\leq  g, \quad |\cq^{n}(\tilde f_\ell)| \leq \alpha^{n} g \quad\text{and}\quad |\cq^{n}(\hat{f}_\ell)| \leq \beta_n \, \alpha^{n} g \quad\text{for all $n,\ell\in \N$,}
\end{equation}
where
\begin{equation*}\label{eq:fhatcrit-S}
\hat{f} = \tilde{f} - \sum_{j \in J} \crr_{j}(f) \quad \text{with} \quad \tilde f= f- \langle \mu, f \rangle.
\end{equation*}

\begin{rem}\label{rem:Guyon} 
In  \cite{Guyon},  only  the   ergodic  Assumptions  \ref{hyp:F}  and \ref{hyp:F1}  were assumed. Indeed, as we will see in Section \ref{sec:special},  if the sequence $\bF = (f_{\ell}, \ell \in \NN)$ is such that $\Pp(f_{\ell}) = 0$ for all $\ell \in \NN,$ then Assumption \ref{hyp:F3} is not needed.  
\end{rem}

\begin{rem}\label{rem:mu-triangle}
We recall that $\mu$ is the invariant probability measure of $\Qq$. We consider the probability measure $\mu^{\vartriangle}$ defined on $S^{3}$ by
\begin{equation*}
\mu^{\vartriangle}(dxx_{0}x_{1}) = \mu(dx)\Pp(x,dx_{0},dx_{1}).
\end{equation*}
Then, for all $f \in \Bb(S^{3})$, one can easily check that $\mu^{\vartriangle}\Qq^{\vartriangle}f = \langle \mu^{\vartriangle}, f \rangle$, that is $\mu^{\vartriangle}$ is the invariant probability measure of $\Qq^{\vartriangle}$. One can also check the following: for all $n\in\NN^{*}$, we have
\begin{equation}\label{eq:Qtrin}
(\Qq^{\vartriangle})^{n}f = \frac{1}{2}(\Qq^{n-1}\Pp f \oplus \Qq^{n-1}\Pp f).
\end{equation}
In the sequel, for all $f \in \Bb(S^{3})$, if $\langle \mu, \Pp f \rangle$ is well defined, we will also set
\begin{equation*}
\tilde{f} = f - \langle \mu^{\vt}, f\rangle = f - \langle \mu, \Pp f \rangle.
\end{equation*}
\end{rem}

The paper is organised as follows. In Section \ref{sec:sub}, we state the main result in the sub-critical case, that is when $2\alpha^{2} < 1$. In Section \ref{sec:special}, we study the special case of conditionally centered functions, that is when $\Pp(f) = 0.$ We will see in particular that in this special case, the value of the ergodicity rate $\alpha \in (0,1)$ does not have any influence on the fluctuations.  In Section \ref{sec:crit-S}, we study the critical and the super-critical cases, that is when $2\alpha^{2} = 1$ and $2\alpha^{2} > 1.$

\section{The sub-critical case: $2\alpha^{2} < 1$.}\label{sec:sub}

Assume that Assumptions \ref{hyp:F} and \ref{hyp:F2} hold. We shall consider sequences $\bF = (f_{\ell}, \ell\in  \N)$ of elements of $\Bb(S^{3})$ such that for all $\ell \in \NN$, $\cp f_{\ell}$, $(\cp f_{\ell}^{2}, \ell \in \NN)$ and $\cp f_{\ell}^{4}$ exist and $(\cp f_{\ell}, \ell\in\NN)$, $(\cp f_{\ell}^{2}, \ell \in \NN)$ and $(\cp f_{\ell}^{4}, \ell \in \NN)$ are sequences of elements of $F$. We shall also assume that \eqref{eq:unif-f} holds for the sequences $(\cp f_{\ell}, \ell\in\NN)$, $(\cp f_{\ell}^{2}, \ell \in \NN)$ and $(\cp f_{\ell}^{4}, \ell \in \NN)$. For a sequence $\bF = (f_{\ell}, \ell\in  \N)$ of elements of $\Bb(S^{3})$, we set 
\begin{equation*}\label{eq:Ntriang}
N_{n,\emptyset}(\bF) = |\GG_{n}|^{-1/2}\sum_{\ell = 0}^{n} M_{\GG_{n-\ell}}(\tilde{f}_{\ell}).
\end{equation*}
where we set $\tilde{f}_{\ell} = f_{\ell} - \langle \mu^{\vartriangle}, f_{\ell}\rangle = f_{\ell} - \langle \mu, \Pp f_{\ell}\rangle.$ 
We have the following result.
\begin{theo}\label{theo:triangsub}
Assume that Assumptions \ref{hyp:F} and \ref{hyp:F2} hold with $\alpha \in (0,1/\sqrt{2})$. For all sequence $\bF = (f_{\ell}, \ell\in  \N)$ of elements of $\Bb(S^{3})$ such that for all $\ell \in \NN$, $\cp f_{\ell}$, $\cp (f_{\ell}^{2})$ and $\cp(f_{\ell}^{4})$ exist, $(\cp f_{\ell}, \ell\in\NN)$, $(\cp f_{\ell}^{2}, \ell \in \NN)$ and $(\cp f_{\ell}^{4}, \ell \in \NN)$ are sequences of elements of $F$ which satisfy \eqref{eq:unif-f} for some $g\in F$, we have the following convergence in distribution
\begin{equation*}
N_{n,\emptyset}(\bF) \; \xrightarrow[n\rightarrow \infty ]{\text{(d)}} \; G,
\end{equation*}
where $G$ is a Gaussian real-valued random variable with
covariance $\stsub (\bF)$ defined by
\[\stsub(\bF) = \stsub_{1}(\bF) + \stsub_{2}(\bF),\]
with $\stsub_{1}(\bF)$ and $\stsub_{2}(\bF)$ are defined by:
\begin{equation*}\label{eq:defS1Sub}
\stsub_{1}(\bF) 
= \sum_{\ell=0}^{+\infty} 2^{-\ell} \langle \mu, \cp ((\tilde{f}_{\ell})^{2})\rangle + \sum_{\ell\geq 0, k\geq 0} 2^{k-\ell} \langle \mu, \cp(\cq^{k}(\cp \tilde{f}_{\ell}) \otimes \cq^{k}(\cp \tilde{f}_{\ell})),
\end{equation*}
and
\begin{equation*}\label{eq:defS2Sub}
\stsub_{2}(\bF) = \stsub_{2,1}(\bF) + \stsub_{2,2}(\bF),
\end{equation*}
with
\begin{equation*}
\stsub_{2,1}(\bF) = \sum_{0\leq \ell < k} 2^{-\ell-1}
\langle \mu, \cp(\tilde{f}_{k}(\cq^{k-\ell - 1}
\cp(\tilde{f}_{\ell}) \oplus \cq^{k-\ell - 1}
\cp(\tilde{f}_{\ell}))) \rangle, 
\end{equation*}
\begin{equation*}
\stsub_{2,2}(\bF) = \sum_{\substack{0\leq \ell< k\\ r\geq
    0}} 2^{r-\ell} \langle \mu,
\cp(\cq^{r}(\cp(\tilde{f}_{k})) \otimes_{sym} \cq^{r +
  k-\ell }(\cp(\tilde{f}_{\ell}))) \rangle. 
\end{equation*}
\end{theo}

\begin{rem}
Let $f \in \Bb(S^{3})$ such that $\Pp f$ and $\Pp(f)^{2}$ exist and belong to $F$. If we take $\bF = (f,0,0,\cdots),$ the infinite sequence where only the first component is non-zero, we obtain
\begin{equation*}
|\GG_{n}|^{-1/2}M_{\GG_{n}}(f) \; \xrightarrow[n\rightarrow \infty ]{\text{(d)}} \; G_{1},
\end{equation*}
where $G_{1}$ is a Gaussian real-valued random variable with covariance $\Sigma^{\vt,sub}_{\GG}(f)$ defined by
\begin{equation*}
\Sigma^{\vt,sub}_{\GG}(f) = \langle \mu,\Pp((\tilde{f})^{2}) \rangle + \sum_{k=0}^{\infty} 2^{k} \langle \mu,\Pp(\Qq^{k}(\Pp\tilde{f} \otimes \Qq^{k}(\Pp\tilde{f})\rangle.
\end{equation*}
If we take $f \in \Bb(S^{3})$ such that $\Pp(f) = 0$ and $\bF = (f, f, \ldots)$, we recover the results of Guyon (see Theorem 19 and Corollary 20 in \cite{Guyon}). In fact, in this special case, we will see in the next section that the fluctuations do not depend on the values of the ergodicity rate $\alpha \in (0,1)$. 
\medskip

Moreover, for two reals $a$ and $b$, if we set $f_{a,b} = a(f - \Pp(f)) + b\Pp\tilde{f}$, we have
\begin{equation*}
|\GG_{n}|^{-1/2}M_{\GG_{n}}(f_{a,b}) \; \xrightarrow[n\rightarrow \infty ]{\text{(d)}} \; G_{1,a,b}
\end{equation*}
where $G_{1,a,b}$ is a Gaussian real-valued random variable with covariance $\Sigma^{\vt,sub}_{\GG}(f_{a,b})$ defined by
\begin{equation*}
\Sigma^{\vt,sub}_{\GG}(f_{a,b}) = a^{2}\langle \mu,\Pp f^{2} - (\Pp f)^{2} \rangle + b^{2} \Sigma_{\GG}^{sub}(\Pp f),
\end{equation*}
where
\begin{equation}\label{eq:SsPft}
\Sigma_{\GG}^{sub}(\Pp f) = \langle \mu,(\Pp \tilde{f})^{2} \rangle + \sum_{k=0}^{\infty} 2^{k} \langle \mu,\Pp(\Qq^{k}(\Pp\tilde{f} \otimes \Qq^{k}(\Pp\tilde{f})\rangle.
\end{equation}
The latter limit implies that
\begin{equation*}
\begin{pmatrix} |\GG_{n}|^{-1/2}M_{\GG_{n}}(f - \Pp f) \\ \\ |\GG_{n}|^{-1/2} M_{\GG_{n}}(\Pp\tilde{f}) \end{pmatrix} \inL G_{2},
\end{equation*}
where $G_{2}$ is a two-dimensional Gaussian vector with zero mean and covariance matrix $\Sigma_{\GG,2}^{\vartriangle}(f)$ defined by
\begin{equation*}
\Sigma_{\GG,2}^{\vartriangle}(f) = \begin{pmatrix} \langle \mu,\Pp f^{2} - (\Pp f)^{2} \rangle &  0 \\ 0  &  \Sigma_{\GG}^{sub}(\Pp f)\end{pmatrix},
\end{equation*}
where $\ssub_{\GG}(\Pp f)$ is defined in \eqref{eq:SsPft}. More generally, for the subtree $\TT_{n}$, we have
\begin{equation*}
\begin{pmatrix} |\TT_{n}|^{-1/2}M_{\TT_{n}}(f - \Pp f) \\ \\ |\TT_{n}|^{-1/2} M_{\TT_{n}}(\Pp\tilde{f}) \end{pmatrix} \inL G_{2,\TT},
\end{equation*}
where $G_{2,\TT}$ is a two-dimensional Gaussian vector with zero mean and covariance matrix $\Sigma_{\TT,2}^{\vartriangle,sub}(f)$ defined by
\begin{equation*}
\stsub_{\TT,2}(f) = \begin{pmatrix} \Sigma_{1,1}(f) & \Sigma_{1,2}(f) \\ \Sigma_{2,1}(f) & \Sigma_{2,2}(f) \end{pmatrix},
\end{equation*}
with
\begin{align*}
&\Sigma_{1,1}(f) = \langle \mu,\Pp((f)^{2} - (\Pp f)^{2}) \rangle,& \\ &\Sigma_{1,2}(f) = \Sigma_{2,1}(f) = \sum_{0\leq \ell < k} 2^{-\ell} \langle \mu,\Pp((f - \Pp f)(\Qq^{k-\ell-1}\Pp\tilde{f} \oplus \Qq^{k-\ell-1}\Pp\tilde{f})) \rangle,& \\ &\Sigma_{2,2}(f) = \Sigma_{\GG}^{sub}(\Pp f) + \sum_{0 \leq \ell < k} 2^{- \ell} \langle \mu,(\Pp\tilde{f})\Qq^{k-\ell}\Pp\tilde{f} \rangle& \\  & \hspace{4cm} + \sum_{\substack{0 \leq \ell < k \\ r \geq 0}} 2^{r-\ell} \langle \mu,\Pp(\Qq^{r}(\Pp\tilde{f}) \otimes_{sym} \Qq^{r+k-\ell}(\Pp\tilde{f})) \rangle,&
\end{align*}
where $\ssub_{\GG}(\Pp f)$ is defined in \eqref{eq:SsPft}.

Now, if we set $f_{a,b} = a(f - \Pp(f)) + b\tilde{f}$ (where we recall that $\tilde{f} = f - \langle \mu,\Pp f \rangle$), then we have
\begin{equation*}
|\GG_{n}|^{-1/2}M_{\GG_{n}}(f_{a,b}) \; \xrightarrow[n\rightarrow \infty ]{\text{(d)}} \; G^{'}_{1,a,b}
\end{equation*}
where $G^{'}_{1,a,b}$ is a Gaussian real-valued random variable with covariance $\Sigma^{\vt,sub'}_{\GG}(f_{a,b})$ defined by
\begin{equation*}
\Sigma^{\vt,sub'}_{\GG}(f_{a,b}) = a^{2}\langle \mu,\Pp((f)^{2} - (\Pp f)^{2}) \rangle + b^{2} \Sigma_{\GG}^{sub}(\Pp f) + 2ab\langle \mu,\Pp((f)^{2} - (\Pp f)^{2}) \rangle,
\end{equation*}
where $\ssub_{\GG}(\Pp f)$ is given in \eqref{eq:SsPft}. The latter limit implies that
\begin{equation*}
\begin{pmatrix} |\GG_{n}|^{-1/2}M_{\GG_{n}}(f - \Pp f) \\ \\ |\GG_{n}|^{-1/2} M_{\GG_{n}}(\tilde{f}) \end{pmatrix} \inL G_{2}^{'},
\end{equation*}
where $G_{2}^{'}$ is a two-dimensional Gaussian vector with zero mean and covariance matrix $\Sigma_{\GG,2}^{\vt'}(f)$ defined by
\begin{equation*}
\Sigma_{\GG,2}^{\vartriangle'}(f) = \begin{pmatrix} \langle \mu,\Pp((f)^{2} - (\Pp f)^{2}) \rangle &  2 \langle \mu,\Pp((f)^{2} - (\Pp f)^{2}) \rangle \\ \\ 2 \langle \mu,\Pp((f)^{2} - (\Pp f)^{2}) \rangle  &  \Sigma_{\GG}^{sub}(\Pp f)\end{pmatrix}.
\end{equation*}
As a consequence, for $f \in \Bb(S^{3})$ such that $\Pp f$ and $\Pp(f^{2})$ exist and belong to $F$, we have that $M_{\GG_{n}}(f - \Pp f)$ and $M_{\GG_{n}}(\Pp f)$ are asymptotically independent, which is not the case for $M_{\TT_{n}}(f - \Pp f)$ and $M_{\TT_{n}}(\Pp f)$, and for $M_{\GG_{n}}(f - \Pp f)$ and $M_{\GG_{n}}(f).$  
\end{rem}
\begin{proof}[Proof of Theorem \ref{theo:triangsub}]
Note that \eqref{eq:Qtrin} and \eqref{eq:geom-erg} imply that for all $f \in \Bb(S^{3})$ such that $\Pp(f) \in F$, there exists $g \in \Bb(S^{3})$ such that
\begin{equation*}
\Pp(g) \in F \quad \text{and} \quad |(\Qq^{\vt})^{n} f - \langle \mu^{\vartriangle},f \rangle| \leq  \alpha^{n} g \quad \text{for all $n \in \NN$}.
\end{equation*}
Applying the results of Theorem 3.1 in Bitseki-Delmas \cite{BD1} to the bifurcating Markov chain $(X_{u}^{\vartriangle}, u \in \TT)$, we get
\[
\stsub_{2}(\bF) = \sum_{0\leq \ell< k} 2^{-\ell} \langle
\mu^{\vartriangle},   \tilde{f}_{k}
(\cq^{\vartriangle})^{k-\ell} \tilde{f}_{\ell}\rangle 
+ \sum_{\substack{0\leq \ell< k\\ 
r\geq 0}} 2^{r-\ell-1} \langle
\mu^{\vartriangle}, \cp^{\vartriangle}\left( (\cq^{\vartriangle})^r
  \tilde{f}_{k} \sot (\cq^{\vartriangle})^{k-\ell+r}
  \tilde{f}_{\ell}\right)\rangle
  \] 
and
\begin{equation*}
\stsub_{1}(\bF) = \sum_{\ell\geq 0} 2^{-\ell} \, \langle
\mu^{\vartriangle},   (\tilde{f}_{\ell})^ 2\rangle +
\sum_{\ell\geq 0, \, k\geq 0}  2^{k-\ell} \, \langle \mu^{\vartriangle},
\cp^{\vartriangle}\left((\cq^{\vartriangle})^k
  \tilde{f}_{\ell} \otimes (\cq^{\vartriangle})^k
  \tilde{f}_{\ell}\right)\rangle. 
\end{equation*}
Using \eqref{eq:Qtrin}, we obtain
\begin{multline*}
\langle \mu^{\vartriangle}, \cp^{\vartriangle}\left( (\cq^{\vartriangle})^r \tilde{f}_{k} \otimes (\cq^{\vartriangle})^{k-\ell+r} \tilde{f}_{\ell}\right)\rangle \\ = \frac{1}{4}\langle \mu^{\vartriangle}, \Pp^{\vartriangle}(\Qq^{r-1}\Pp\tilde{f}_{k} \oplus \Qq^{r-1}\Pp\tilde{f}_{k}) \otimes (\Qq^{k - \ell + r-1}\Pp\tilde{f}_{\ell} \oplus \Qq^{k - \ell + r-1}\Pp\tilde{f}_{\ell})\rangle.
\end{multline*}
Next, using the definition of $\Pp^{\vartriangle}$, we obtain
\begin{multline*}
\frac{1}{4}(\Pp^{\vartriangle}(\Qq^{r-1}\Pp\tilde{f}_{k} \oplus \Qq^{r-1}\Pp\tilde{f}_{k}) \otimes (\Qq^{k - \ell + r-1}\Pp\tilde{f}_{\ell} \oplus \Qq^{k - \ell + r-1}\Pp\tilde{f}_{\ell})) \\ = \frac{1}{4}(2\Qq(\Qq^{r-1}\Pp(\tilde{f}_{k})) \otimes 2\Qq(\Qq^{k - \ell + r-1}\Pp(\tilde{f}_{k}))) = \Qq^{r}(\Pp\tilde{f}_{k}) \otimes \Qq^{k - \ell + r} \Pp(\tilde{f}_{\ell}).
\end{multline*}
From the two previous equalities and the definition of $\mu^{\vartriangle}$ we are led to
\begin{align}
\langle \mu^{\vartriangle}, \cp^{\vartriangle}\left( (\cq^{\vartriangle})^r \tilde{f}_{k} \otimes (\cq^{\vartriangle})^{k-\ell+r} \tilde{f}_{\ell}  \right)\rangle &= \langle \mu^{\vartriangle}, \Qq^{r}(\Pp\tilde{f}_{k}) \otimes \Qq^{k - \ell + r} \Pp(\tilde{f}_{\ell})\rangle& \nonumber \\ &=  \langle \mu, \Pp(\Qq^{r}(\Pp\tilde{f}_{k}) \otimes \Qq^{k - \ell + r} \Pp(\tilde{f}_{\ell}))\rangle.& \label{eq:defS22tri}
\end{align}
In the same way, from \eqref{eq:Qtrin} and using the definition of $\mu^{\vartriangle}$, we are led to
\begin{align}
\langle \mu^{\vartriangle},   \tilde{f}_{k} (\cq^{\vartriangle})^{k-\ell} \tilde{f}_{\ell}\rangle &= \langle \mu^{\vartriangle},\tilde{f}_{k}((1/2)(\Qq^{k-\ell-1}\Pp(\tilde{f}_{\ell}) \oplus \Qq^{k-\ell-1}\Pp(\tilde{f}_{\ell})) \rangle& \nonumber \\ &= (1/2) \langle \mu,\Pp(\tilde{f}_{k}(\Qq^{k-\ell-1}\Pp(\tilde{f}_{\ell}) \oplus \Qq^{k-\ell-1}\Pp(\tilde{f}_{\ell}))) \rangle.& \label{eq:defS21tri}
\end{align}
Now, from \eqref{eq:defS22tri} and \eqref{eq:defS21tri} we obtain the
definition of $\stsub_{2}(\bF)$ given in Theorem
\ref{theo:triangsub}. We can obtain
$\stsub_{1}(\bF)$ in the same way and this ends
the proof. 
\end{proof}

\section{The special case of sequences $(f_{\ell}, \ell \in \NN)$ such that $\Pp f_{\ell} = 0$ for all $\ell \in \NN.$}\label{sec:special}

Under the assumptions of Theorem \ref{theo:triangsub}, let $\bF = (f_{\ell},\ell\in\NN)$ be a sequence such that $\Pp f_{\ell} = 0$ for all $\ell \in \NN$. We will show that in this special case, the fluctuations does not depend to the value of the ergodicity rate $\alpha \in (0,1)$. For that purpose we have the following result.
\begin{theo}\label{theo:triangPf}
Under the assumptions of Theorem \ref{theo:triangsub}, let $\bF = (f_{\ell},\ell\in\NN)$ be a sequence such that $\Pp f_{\ell} = 0$ for all $\ell \in \NN$. Then for all $\alpha\in(0,1)$, we have the following convergence in distribution:
\begin{equation*}
N_{n,\emptyset}(\bF) \; \xrightarrow[n\rightarrow \infty ]{\text{(d)}} \; G,
\end{equation*}
where $G$ is a Gaussian real-valued random variable with covariance $\Sigma(\bF)$ defined by \eqref{eq:cv-Vn-Pf0}.
\end{theo}
\begin{rem}\label{rem:ext-GDM}
We stress that Theorem \ref{theo:triangPf} allows to recover the results of Guyon \cite{Guyon} and the results of Delmas and Marsalle \cite{DM10}. In fact, on the one hand, if we take $\bF = (f,f,\cdots)$, the infinite sequence of the same function $f$, we get 
\begin{equation*}
|\GG_{n}|^{-1/2} M_{\TT_{n}}(f) \inL G,
\end{equation*}
where $G$ is the Gaussian real-valued random variable with mean $0$ and variance $\Sigma_{\TT}^{\vt} = 2\langle \mu, \cp f^{2}\rangle$. In particular, using the fact that $\lim_{n\rightarrow\infty}(\TT_{n}/\GG_{n}) = 2$, the latter result implies that
\begin{equation*}\label{eq:cltG}
|\TT_{n}|^{-1/2}M_{\TT_{n}}^{\vt}(f - \cp f) \inL \langle \mu,\cp f^{2} \rangle
\end{equation*}
and we thus recover the results in \cite{Guyon}. Next, let $d\in \NN^{*}$. If in Theorem \ref{theo:triangPf} we take $\bF = (f_{1},\cdots,f_{d},0,0,\cdots)$ where $(f_{1},\cdots,f_{d})\in\Bb(S^{3})^{d}$, with $\cp f_{\ell} = 0$ and $\Pp(f_{\ell}^{2}) \in F$ for all $\ell \in \{1,\cdots,d\}$, we get
\begin{equation*}
|\GG_{n}|^{-1/2} \sum_{\ell = 0}^{d-1} M_{\GG_{n-\ell}}(f_{\ell+1}) \inL G, 
\end{equation*}
where $G$ is a Gaussian real-valued random variable with mean $0$ and variance $$\Sigma_{\GG}^{\vt} = \sum_{\ell = 0}^{d-1} 2^{-\ell}\langle \mu, \cp((f_{\ell+1})^{2})\rangle.$$ 
We thus recover, as is \cite{DM10}, that the fluctuations over different generations are asymptotically independent.
\end{rem}

\begin{proof}
Let $(p_{n}, n\in \NN^{*})$ be an increasing sequence of elements of $\NN^{*}$ which satisfies for all $\lambda>0$:
\begin{equation}\label{eq:def-pn}
p_n< n, \quad \lim_{n\rightarrow \infty } p_n/n=1 \quad\text{and}\quad \lim_{n\rightarrow \infty } n-p_n - \lambda \log(n) = +\infty .
\end{equation}
We have
\[N_{n,\emptyset}(\bF) = R_{0}(n) + \Delta_{n}(\bF),\]
where
\begin{equation*}
R_{0}(n) = |\GG_{n}|^{-1/2}\sum_{\ell = 0}^{n-p-1} \sum_{u\in\GG_{\ell}}f_{n-\ell}(X^{\vartriangle}_{u}) \quad \text{and} \quad \Delta_{n}(\bF) = \sum_{i\in\GG_{n-p}} \Delta_{n,i}(\bF),
\end{equation*}
with the $\cf_{i}\text{-}$martingale increments $\Delta_{n,i}(\bF)$ defined by $$\Delta_{n,i}(\bF) =  |\GG_{n}|^{-1/2}\sum_{\ell=0}^{p} \sum_{u\in i\GG_{p-\ell}} f_{\ell}(X_{u}^{\vartriangle}).$$ 

We have
\begin{align}
\EE[R_{0}(n)^{2}]  &= \frac{1}{|\GG_{n}|}\EE[(\sum_{\ell = 0}^{n-p-1} \sum_{u\in \GG_{\ell}} \tilde{f}_{n-\ell}(X^{\vt}_{u}))^{2}]& \nonumber\\ &\leq \frac{1}{|\GG_{n}|} (\sum_{\ell = 0}^{n-p-1} \EE[(\sum_{u\in\GG_{\ell}}\tilde{f}_{n-\ell}(X^{\vt}_{u}))^{2}]^{1/2})^{2},& \label{eq:BR0Triang}
\end{align}
where we use the Minkowski inequality for the first inequality. By developing the term in the expectation, we get
\begin{align*}
\EE[(\sum_{u\in\GG_{\ell}}\tilde{f}_{n-\ell}(X^{\vt}_{u}))^{2}] &= \EE[\sum_{u\neq v\in\GG_{\ell}} \EE[\tilde{f}_{n-\ell}(X^{\vt}_{u}) \tilde{f}_{n-\ell}(X^{\vt}_{v})|X_{u},X_{v}]] + \EE[\sum_{u\in\GG_{\ell}} \EE[\tilde{f}_{n-\ell}^{2}(X^{\vt}_{u})|X_{u}]]& \\ &= \EE[\sum_{u\neq v\in\GG_{\ell}} \cp\tilde{f}_{n-\ell}(X_{u}) \cp\tilde{f}_{n-\ell}(X_{v})] + \EE[\sum_{u\in\GG_{\ell}} \cp(\tilde{f}_{n-\ell}^{2})(X_{u})]& \\ & = \EE[(\sum_{u\in\GG_{\ell}}\cp\tilde{f}_{n-\ell}(X_{u}))^{2}] + \EE[\sum_{u\in\GG_{\ell}} (\cp(f_{n-\ell}^{2}) - (\cp f_{n-\ell})^{2})(X_{u})],& 
\end{align*}
where we used the branching Markov property for the second inequality and the fact that $\cp(\tilde{f}_{n-\ell}^{2}) - (\cp \tilde{f}_{n-\ell})^{2} = \cp(f_{n-\ell}^{2}) - (\cp f_{n-\ell})^{2}$ for the third equality. Using this last equality in \eqref{eq:BR0Triang} and using the inequalities $\sqrt{a+b} \leq \sqrt{a}+\sqrt{b}$ and $(a+b)^{2} \leq 2a^{2} + 2b^{2}$, we get
\begin{multline*}\label{eq:bdRtriSub}
\EE[R_{0}(n)^{2}] \leq \frac{1}{|\GG_{n}|}(\sum_{\ell = 0}^{n-p-1} (\EE[(\sum_{u\in\GG_{\ell}}\cp\tilde{f}_{n-\ell}(X_{u}))^{2}]^{1/2} + \EE[\sum_{u\in\GG_{\ell}} (\cp(f_{n-\ell}^{2}) - (\cp f_{n-\ell})^{2})(X_{u})]^{1/2}))^{2} \\ \leq \frac{2}{|\GG_{n}|}(\sum_{\ell = 0}^{n-p-1} \EE[M_{\GG_{\ell}}(\cp(\tilde{f}_{n-\ell}))^{2}]^{1/2})^{2} + \frac{2}{|\GG_{n}|}(\sum_{\ell = 0}^{n-p-1}\EE[M_{\GG_{\ell}}(\cp(f_{n-\ell}^{2}) - (\cp f_{n-\ell})^{2})]^{1/2})^{2}.
\end{multline*}
On the one hand, by replacing $\tilde{f}_{n-\ell}$ with $\cp(\tilde{f}_{n-\ell})$ in the Lemma \ref{lem:cvR0-M}, we get
\begin{equation}\label{eq:lR0tri1}
\lim_{n\rightarrow +\infty}\frac{2}{|\GG_{n}|}(\sum_{\ell = 0}^{n-p-1} \EE[M_{\GG_{\ell}}(\cp(\tilde{f}_{n-\ell}))^{2}]^{1/2})^{2} = 0.
\end{equation}
On the other hand, from \eqref{eq:Q1}, \eqref{eq:unif-f} and $(iii)$ of Assumption \ref{hyp:F}, we get that there is a positive constant $c$ independent of $n$ such that
\begin{equation*}
 \EE[M_{\GG_{\ell}}(\cp(f_{n-\ell}^{2}) - (\cp f_{n-\ell})^{2})] \leq c2^{\ell}.
\end{equation*}
The latter inequality implies that
\begin{align*}
\frac{2}{|\GG_{n}|}(\sum_{\ell = 0}^{n-p-1}\EE[M_{\GG_{\ell}}(\cp(f_{n-\ell}^{2}) - (\cp f_{n-\ell})^{2})]^{1/2})^{2} &\leq \frac{2}{|\GG_{n}|}(\sum_{\ell=0}^{n-p-1}c^{1/2}2^{\ell/2})& \\ &\leq \frac{2c}{(\sqrt{2} - 1)^{2}} 2^{-p}.&
\end{align*}
Since $\lim_{n\rightarrow +\infty} p_{n} = +\infty$, we deduce that
\begin{equation}\label{eq:lR0tri2}
\lim_{n\rightarrow +\infty}\frac{2}{|\GG_{n}|}(\sum_{\ell = 0}^{n-p-1}\EE[M_{\GG_{\ell}}(\cp(f_{n-\ell}^{2}) - (\cp f_{n-\ell})^{2})]^{1/2})^{2} = 0.
\end{equation}
Finally, from \eqref{eq:lR0tri1} and \eqref{eq:lR0tri2}, we get that
\[\lim_{n\rightarrow \infty} \EE[(R_{0}(n))^{2}] = 0 \quad \text{for all $\alpha \in (0,1)$}.\]

Now, we consider the bracket
\begin{equation*} 
V(n) = \sum_{i \in \GG_{n-p}} \EE[(\Delta_{n,i}(\bF))^{2}|\Ff_{i}].
\end{equation*}
For all $i \in \GG_{n-p}$, we have
\begin{equation*}
(\Delta_{n,i}(\bF))^{2} = |\GG_{n}|^{-1} \sum_{\ell = 0}^{p} M_{i\GG_{p-\ell}}(f_{\ell})^{2} + 2 |\GG_{n}|^{-1} \sum_{0 \leq \ell < k \leq p} M_{i\GG_{p-\ell}}(f_{\ell}) M_{i\GG_{p-k}}(f_{k}).
\end{equation*}
Using the branching Markov chain property, we have, for all $0\leq \ell < k \leq p$:
\begin{equation*}
\EE[M_{i\GG_{p-\ell}}(f_{\ell}) M_{i\GG_{p-k}}(f_{k})|\Ff_{i}] = \EE_{X_{i}}[M_{\GG_{p-k}}(f_{k}) M_{\GG_{p-\ell}}(\Pp f_{\ell})] = 0,
\end{equation*}
where we used the fact that $\Pp(f_{\ell}) = 0.$ On the other hand, using again the branching Markov property twice and \eqref{eq:Q1}, we have for all $0\leq \ell \leq p,$
\begin{equation*}
\EE[M_{i\GG_{p-\ell}}(f_{\ell})^{2}|\Ff_{i}] = \EE_{X_{i}}[M_{\GG_{p - \ell}}(f_{\ell}^{2})] = 2^{p-\ell} \Qq^{p-\ell} \Pp f_{\ell}^{2}(X_{i}).
\end{equation*}
From the two latter equalities, it follows that
\begin{equation*}
V_{n} = |\GG_{n-p}|^{-1}\sum_{i \in \GG_{n-p}} \sum_{\ell = 0}^{p} 2^{-\ell} \Qq^{p-\ell} \Pp f_{\ell}^{2}(X_{i}).
\end{equation*}
We set
\begin{equation*}
\Sigma(\bF) = \sum_{\ell = 0}^{\infty} 2^{-\ell} \langle \mu, \Pp f^{2}_{\ell} \rangle \quad \text{and} \quad \Sigma_{n}(\bF) = \sum_{\ell = 0}^{p} 2^{-\ell} \langle \mu, \Pp f^{2}_{\ell} \rangle.
\end{equation*}
On the one hand, we have
\begin{align}
\EE[|V_{n} - \Sigma_{n}(\bF)|] &\leq \sum_{\ell = 0}^{p} 2^{-\ell} |\GG_{n-p}|^{-1} \EE[M_{\GG_{n-p}}(|\Qq^{p-\ell}(\Pp f_{\ell}^{2} - \langle \mu, \Pp f_{\ell}^{2} \rangle|)] \nonumber \\ & \leq \sum_{\ell = 0}^{p} 2^{-\ell} \alpha^{p-\ell} \langle \nu,g\rangle \nonumber \\ & \leq c \sum_{\ell = 0}^{p} 2^{-\ell} \alpha^{p-\ell}, \label{eq:IVn-SnPf0}
\end{align}
where we used \eqref{eq:geom-erg} for the second inequality and $(iii)$ of Assumption \ref{hyp:F} for the last inequality. Now, since
\begin{equation*}
\sum_{\ell = 0}^{\infty} 2^{-\ell} < \infty \quad \text{and} \quad \lim_{n \rightarrow \infty} \alpha^{p-\ell} = 0 \quad \text{for all fixed $\ell$},
\end{equation*}
we get, using the dominated convergence theorem and \eqref{eq:IVn-SnPf0},
\begin{equation*}
\lim_{n\rightarrow \infty} \EE[|V_{n} - \Sigma_{n}(\bF)|] = 0.
\end{equation*}
On the other hand, since $(\Pp f_{\ell}^{2}, \ell \in \NN)$ satisfies uniformly \eqref{eq:erg-bd} and \eqref{eq:geom-erg}, it follows that
\begin{equation*}
\lim_{n \rightarrow \infty} |\Sigma(\bF) - \Sigma_{n}(\bF)| = 0.
\end{equation*}
From the foregoing, we deduce that, for all $\alpha \in (0,1)$, 
\begin{equation}\label{eq:cv-Vn-Pf0}
\lim_{n \rightarrow \infty} V_{n} = \Sigma(\bF) \quad \text{in probability}.
\end{equation}

Finally, in order to check the Lindeberg condition, we consider $R_{3}(n) = \sum_{i\in \GG_{n-p}} \EE[\Delta_{n,i}(\bF)^{4}]$. Using that $(\sum_{k=0}^p a_k)^4 \leq  (p+1)^3 \sum_{k=0}^p a_k^4$, we get
\begin{equation}\label{eq:I-R3-Pf0}
R_{3}(n) \leq (p+1)^{3} |\GG_{n}|^{-2} \sum_{\ell = 0}^{p} \sum_{i \in \GG_{n-p}} \EE[(M_{i\GG_{p-\ell}}(f_{\ell}))^{4}].
\end{equation}
For all $i \in \GG_{n-p}$ and for all $0\leq \ell \leq p$, we have, using the branching Markov property (see the calculus in \cite{BDG14}, Remark 2.3 for more details), \eqref{eq:Q1}, \eqref{eq:Q2} and \eqref{eq:erg-bd}:
\begin{align*}
\EE[M_{i\GG_{p-\ell}}(f_{\ell})^{4}] & \leq \EE[M_{i\GG_{p-\ell}}(\Pp f_{\ell}^{4})] + 6 \, \EE[M_{i\GG_{p - \ell}}(\Pp f_{\ell}^{2})^{2}] \\ & = 2^{p-\ell} \, \EE[\Qq^{p-\ell}\Pp f_{\ell}^{4}(X_{i})] + 6 \, \EE[2^{p-\ell}\Qq^{p-\ell}(\Pp f_{\ell}^{2})^{2}(X_{i})]  \\ & \hspace{0.75cm} + 6 \,\EE[\sum_{k = 0}^{p - \ell - 1} 2^{p-\ell+k}\Qq^{p-\ell-k-1}(\Pp(\Qq^{k}\Pp f_{\ell}^{2} \otimes^{2}))(X_{i})] \\ & \leq c \, 2^{p-\ell} \, \E[(g_{1}(X_{i}) + 2^{p-\ell}g_{2}(X_{i}))],
\end{align*}
where $g_{1},g_{2} \in F.$ Using \eqref{eq:I-R3-Pf0}, the latter inequality, \eqref{eq:Q1} and $(iii)$ of Assumption \ref{hyp:F}, we get
\begin{align*}
R_{3}(n) &\leq c \, p^{3}|\GG_{n}|^{-2} \sum_{\ell = 0}^{p} 2^{p-\ell} \EE[M_{\GG_{n-p}}(g_{1})] + c \, p^{3}|\GG_{n}|^{-2} \sum_{\ell = 0}^{p} 2^{2(p-\ell)} \EE[M_{\GG_{n-p}}(g_{2})] \\ & \leq c \, p^{3} \, 2^{-n} + c \, p^{3} \, 2^{- n + p}.
\end{align*}
Using \eqref{eq:def-pn}, it follows from the foregoing that for all $\alpha \in (0,1)$, 
\[\lim_{n \rightarrow \infty} R_{3}(n) = 0.\]
\end{proof}
 
\section{The critical and the super-critical cases: $2\alpha^{2} = 1$ and $2\alpha^{2} > 1$.}\label{sec:crit-S}

First note that
\begin{equation*}
N_{n,\emptyset}(\bF) = |\GG_{n}|^{-1/2} \sum_{\ell = 0}^{n} M_{\GG_{n-\ell}}(f_{\ell} - \Pp f_{\ell}) + |\GG_{n}|^{-1/2} \sum_{\ell = 0}^{n} M_{\GG_{n-\ell}}(\Pp f_{\ell} - \langle \mu, \Pp f_{\ell}\rangle).
\end{equation*}
Next, let $(s_{n}, n\in \NN)$ be a sequence of real numbers converging to 0. From Theorem \ref{theo:triangPf}, we have
\begin{equation}\label{eq:r-cri-scri}
s_{n}|\GG_{n}|^{-1/2} \sum_{\ell = 0}^{n} M_{\GG_{n-\ell}}(f_{\ell} - \Pp f_{\ell}) \nlim 0 \quad \text{in probability},
\end{equation}
so that in the critical and the super-critical case, the study of the asymptotic law of $N_{n,\emptyset}(\bF)$ is reduced to the study of the asymptotic law of $|\GG_{n}|^{-1/2} \sum_{\ell = 0}^{n} M_{\GG_{n-\ell}}(\Pp f_{\ell} - \langle \mu, \Pp f_{\ell}\rangle).$

We introduce the following notation for $k,\ell \in \NN$:
\begin{equation*}\label{eq:Pfkltri}
\cp f_{k,\ell}^{*} = \sum_{j\in J} \theta_{j}^{\ell - k} \crr_{j}(\cp f_{k}) \sot \overline{\crr}_{j}(\cp f_{\ell}).
\end{equation*}

Then, from \eqref{eq:r-cri-scri} and Theorem 3.4 in Bitseki-Delmas \cite{BD1} applied to the functions $(\Pp f_{\ell}, \ell \in \NN)$, we get the following.
\begin{theo}\label{theo:triangcrit}
Assume that Assumptions \ref{hyp:F} and \ref{hyp:F3} hold with $\alpha = 1/\sqrt{2}$. For all sequences $\bF = (f_{\ell}, \ell \in  \NN)$ of elements of $\Bb(S^{3})$ such that  for all $\ell \in \NN$, $\cp f_{\ell}$ and $\cp f_{\ell}^{2}$ exist, $(\cp f_{\ell}, \ell\in\NN)$, $(\cp f_{\ell}^{2}, \ell \in \NN)$ and $(\cp f_{\ell}^{4}, \ell \in \NN)$ are sequences of elements of $F$ which satisfy \eqref{eq:unif-f-crit} for some $g\in F$, we have the following convergence in distribution:
\begin{equation*}
n^{-1/2}N_{n,\emptyset}(\bF) \; \xrightarrow[n\rightarrow \infty ]{\text{(d)}} \; G,
\end{equation*}
where $G$ is a Gaussian real-valued random variable with covariance $\stcrit(\bF^{\vartriangle})$ defined by
\begin{equation*}
\stcrit(\bF) = \stcrit_{1}(\bF) + \stcrit_{2}(\bF),
\end{equation*}
where 
\begin{equation*}\label{eq:defS1crit}
\stcrit_{1}(\bF) = \sum_{\ell = 0}^{+\infty}
2^{-\ell}\langle \mu,\cp(\cp f_{\ell,\ell}^{*})\rangle
\quad \text{and} \quad 
\stcrit_{2}(\bF^{\vartriangle}) = \sum_{0\leq \ell < k} 2^{-(k+\ell)/2}\langle \mu,\cp(\cp f_{k,\ell}^{*})\rangle.
\end{equation*}
\end{theo}
\begin{rem}\label{rem:theotricrit}
We stress that the covariance $\Sigma^{(\vartriangle,crit)}(\bF)$ can take the value $0$. This is the case if for all $\ell \in \NN$, the orthogonal projection of $\cp(f_{\ell})$ on the eigen-spaces corresponding to the eigenvalues $1$ and $\alpha_{j}$, $j \in J$, equal 0. However, the case where $\cp(f_{\ell}) = 0$ for all $\ell \in \NN$ is treated in Theorem \ref{theo:triangPf}. Indeed, we have seen that in this case the good normalisation is not $(n|\GG_{n}|)^{-1/2}$, but $|\GG_{n}|^{-1/2}$.
\end{rem}

Next, from \eqref{eq:r-cri-scri} and Theorem 3.9 in Bitseki-Delmas \cite{BD1} applied to the functions $(\Pp f_{\ell}, \ell \in \NN)$, we get the following result. 
\begin{theo}\label{theo:triScrit}
Assume that Assumptions \ref{hyp:F} and \ref{hyp:F3} hold with $\alpha > 1/\sqrt{2}$. For all sequences $\bF = (f_{\ell}, \ell \in  \N)$ of elements of $\Bb(S^{3})$ such that  for all $\ell \in \NN$, $\cp f_{\ell}$ and $\cp f_{\ell}^{2}$ exist, $(\cp f_{\ell}, \ell \in \NN)$ and $(\cp f_{\ell}^{2}, \ell \in \NN)$ are sequences of elements of $F$ which satisfy \eqref{eq:unif-f-crit} for some $g\in F$, we have the following convergence in probability:
\begin{equation*}
(2\alpha^{2})^{-n/2} N_{n, \emptyset}(\bF)
- \sum_{\ell\in \N} (2\alpha)^{-\ell} \sum_{j\in J} \theta_j^{n-\ell}
M_{\infty  ,j}(\cp(f_{\ell}))  
\; \xrightarrow[n\rightarrow \infty ]{\P}  \; 0. 
\end{equation*}
\end{theo}
\begin{rem}\label{rem:theo-triScrit}
If the sequence $\bF = (f_{\ell},\ell \in \NN)$ is such that for all $\ell \in \NN$, the orthogonal projection of $\cp(f_{\ell})$ on the eigen-spaces corresponding to the eigenvalues $1$ and $\alpha_{j}$, $j \in J$, equal 0, then we have the following convergence in probability:
\begin{equation*}
(2\alpha^{2})^{-n/2} N_{n, \emptyset}(\bF)  \; \xrightarrow[n\rightarrow \infty ]{\P}  \; 0. 
\end{equation*}
Once again, the case where $\cp(f_{\ell}) = 0$ for all $\ell \in \NN$ is treated in Theorem \ref{theo:triangPf}. Indeed, we have seen that in this case the good normalisation is not $((2\alpha^{2})^{n}|\GG_{n}|)^{-1/2}$, but $|\GG_{n}|^{-1/2}$.
\end{rem}

\section{Appendix}

\subsection*{Moments formula for BMC}\label{sec:moment}
Let $X=(X_i,  i\in \T)$ be a  BMC on $(S, \cs)$  with probability kernel $\cp$. Recall that  $|\G_n|=2^n$   and $M_{\G_n}(f)=\sum_{i\in \G_n} f(X_i)$. We also recall that $2\cq(x,A)=\cp(x, A\times S) + \cp(x, S\times A)$ for $A\in \cs$. We use the convention that $\sum_\emptyset=0$.

\medskip

We recall the following well known and easy to establish  many-to-one formulas  for BMC. 

\begin{lem}\label{lem:Qi}
Let $f,g\in \cb(S)$, $x\in S$ and $n\geq m\geq 0$. Assuming that all the quantities below are well defined,  we have:
\begin{align}\label{eq:Q1}
\E_x\left[M_{\G_n}(f)\right] &=|\G_n|\, \cq^n f(x)= 2^n\, \cq^n f(x) ,\\\label{eq:Q2} \E_x\left[M_{\G_n}(f)^2\right] &=2^n\, \cq^n (f^2) (x) + \sum_{k=0}^{n-1} 2^{n+k}\,   \cq^{n-k-1}\left( \cp \left(\cq^{k}f\otimes \cq^k f \right)\right) (x).
\end{align}
\end{lem}

We have the following lemma.

\begin{lem} \label{lem:cvR0-M}
Under the assumptions of Theorem \ref{theo:triangsub}, let $\bF = (f_{\ell},\ell\in\NN)$ be a sequence of elements of $F$ satisfying Assumptions \ref{hyp:F2} uniformly, that is \reff{eq:unif-f} for some $g\in F$. Let $(p_{n}, n\in \NN^{*})$ be an increasing sequence of elements of $\NN^{*}$ which satisfies \eqref{eq:def-pn}. Then, for all $\alpha \in (0,1)$,
\begin{equation*}
\lim_{n \rightarrow \infty} |\G_n|^{-1/2}\, \sum_{k=0}^{n-p-1} \E[M_{\G_{k}} (\tilde f_{n-k})^2]^{1/2} = 0
\end{equation*}
\end{lem}

\begin{proof}
For all $k\geq 1$, we have:
\begin{align*}
\E_x[M_{\G_{k}} (\tilde f_{n-k})^2] & \leq  2^{k}\,  g_{1}(x)  + \sum_{\ell=0}^{k-1}  2^{k+\ell}\alpha^{2\ell}\, \, \cq^{k-\ell -1} \left(\cp ( g_{2}\otimes g_{2})\right)(x)\\ & \leq  2^{k} \,  g_{1}(x) +   2^{k} \sum_{\ell=0}^{k-1} (2\alpha^2)^{\ell}  g_{3}(x)\\ & \leq 2^{k}g_{4}(x) \times \begin{cases} 1 & \text{if $2\alpha^{2} < 1$} \\ k & \text{if $2\alpha^{2} = 1$} \\ (2\alpha^{2})^{k} & \text{if $2\alpha^{2} > 1,$} \end{cases}
\end{align*}
with $g_1, g_2, g_3, g_4\in F$ and where we used \reff{eq:Q2}, \reff{eq:unif-f} twice and \reff{eq:erg-bd} twice (with $f$ and $g$ replaced by $2(g^2+ \langle  \mu, g \rangle^2)$ and $g_{1}$, and with  $f$ and $g$ replaced by $g$ and $g_{2}$)  for the first inequality, \reff{eq:erg-bd} (with $f$ and $g$  replaced by $\cp ( g_{2}\otimes  g_{2})$ and $g_{3}$) for the second, and $g_4 = g_1+ (1 + |1-2\alpha^2|^{-1}) g_3$ for the last. We deduce that:
\begin{align*}\label{eq:majo-R0}
|\G_n|^{-1/2}\, \sum_{k=0}^{n-p-1} \E[M_{\G_{k}} (\tilde f_{n-k})^2]^{1/2} & \leq c\, 2^{-n/2} \times \begin{cases} \sum_{k=0}^{n-p-1} 2^{k/2} & \text{if $2\alpha^{2} < 1$} \\ \sum_{k=0}^{n-p-1} k^{1/2}2^{k/2} & \text{if $2\alpha^{2} = 1$} \\ \sum_{k=0}^{n-p-1} (2\alpha^{2})^{k/2}2^{k/2} & \text{if $2\alpha^{2} < 1$}\end{cases} \\ &\leq  3c\, \begin{cases} 2^{- p/2} & \text{if $2\alpha^{2} = 1$} \\ (n-p)^{1/2} \, 2^{- p/2} & \text{if $2\alpha^{2} = 1$} \\ (2\alpha^{2})^{(n-p)/2} \, 2^{- p/2} & \text{if $2\alpha^{2} > 1$}. \end{cases}
\end{align*}
Use that $\lim_{n\rightarrow \infty } p_{n}=\infty$ and $\lim_{n\rightarrow \infty } p_n/n=1$ to conclude. 
\end{proof}

\bibliographystyle{abbrv}
\bibliography{biblio}
\end{document}